\documentclass{llncs}
\usepackage{makeidx}  
\usepackage{amsmath}
\usepackage{amssymb}
\usepackage{amsfonts}
\usepackage[english]{babel}
\usepackage{enumerate}
\usepackage{fancyhdr}

\def \egal {\stackrel{{\rm def}}{=}}

\def \egal {\stackrel{{\rm def}}{=}}

\newcommand \cA{{\cal A}}
\newcommand \cB{{\cal B}}
\newcommand \cC{{\cal C}}

\newcommand \cF{{\cal F}}

\newcommand \cP{{\cal P}}

\newcommand \cS{{\cal S}}
\newcommand \cT{{\cal T}}

\newcommand \cX{{\cal X}}
\newcommand \cY{{\cal Y}}
\newcommand \cZ{{\cal Z}}
\newcommand{\1}{{\rm 1}\kern-0.24em{\rm I}}

\newcommand{\DS}{\displaystyle}

\newtheorem{rem}{Remark}

\begin{document}

\title{Suboptimality of Penalized Empirical Risk Minimization in Classification.}
\titlerunning{Optimality of Aggregation Method in Classification}  
%
\author{Guillaume Lecu\'e \footnote{Paper to be considered for the Mark Fulk Award for the "best student paper".}}
\authorrunning{Guillaume Lecu\'e}   
%
\tocauthor{Guillaume Lecu\'e (Universit\'e Paris VI)}
\institute{Laboratoire de Probabilit\'es et Mod\`eles Al\'eatoires
(UMR CNRS 7599)\\Universit\'e Paris VI\\ 4 pl.Jussieu, BP 188,
75252 Paris, France ,\\
\email{lecue@ccr.jussieu.fr}}

\maketitle              

\begin{abstract} \noindent
Let $\cF$ be a set of $M$ classification procedures with values in
$[-1,1]$. Given a loss function, we want to construct a procedure
which mimics at the best possible rate the best procedure in $\cF$.
This fastest rate is called optimal rate of aggregation. Considering
a continuous scale of loss functions with various types of
convexity, we prove that optimal rates of aggregation can be either
$((\log M)/n)^{1/2}$ or $(\log M)/n$. We prove that, if all the $M$
classifiers are binary, the (penalized) Empirical Risk Minimization
procedures are suboptimal (even under the margin/low noise
condition) when the loss function is somewhat more than convex,
whereas, in that case, aggregation procedures with exponential
weights achieve the optimal rate of aggregation.
\end{abstract}

\section{Introduction}\label{SubERMSectionIntroduction}
Consider the problem of binary classification. Let $(\cX,\cA)$ be a
measurable space. Let $(X,Y)$ be a couple of random variables, where
$X$ takes its values in $\cX$ and $Y$ is a random label taking
values in $\{-1,1\}$. We denote by $\pi$ the probability
distribution of $(X,Y)$. For any function
$\phi:\mathbb{R}\longmapsto\mathbb{R},$ define the $\phi-$risk of a
real valued classifier $f:\cX\longmapsto\mathbb{R}$ by
$$A^{\phi}(f)=\mathbb{E}[\phi(Yf(X))].$$ Many different losses
have been discussed in the literature along the last decade (cf.
\cite{cv:95,fs:97,lv:04,fht:00,by:02}), for instance:
\begin{equation*}
  \begin{array}{ll}
  \phi_0(x)=\1_{(x\leq0)} & \mbox{classical loss or }0-1 \mbox{ loss}\\
  \phi_1(x)=\max(0,1-x) & \mbox{hinge loss (SVM loss)}\\
  x\longmapsto\log_2(1+\exp(-x)) & \mbox{logit-boosting loss}\\
  x\longmapsto\exp(-x) & \mbox{exponential boosting loss}\\
  x\longmapsto(1-x)^2 & \mbox{squared loss}\\
  x\longmapsto\max(0,1-x)^2 & \mbox{$2$-norm soft margin loss}\\
  \end{array}
\end{equation*}
We will be especially interested in losses having convex properties
as it is considered in the following definition (cf. \cite{jrt:06}).
\begin{definition}
  Let $\phi:\mathbb{R}\longmapsto\mathbb{R}$ be a function and $\beta$ be a positive
  number. We say that $\phi$ is ${\mathbf{\beta-}}${\bf{convex on}} $[-1,1]$
  when $$[\phi'(x)]^2\leq \beta \phi''(x),\quad \forall |x|\leq 1.$$
\end{definition} For example, logit-boosting loss is $(e/\log 2)-$convex,
exponential boosting loss is $e-$convex, squared and $2-$norm soft
margin losses are $2-$convex.

We denote by $f^{*}_\phi$ a function from $\cX$ to $\mathbb{R}$
which minimizes $A^\phi$ over all real-valued functions and by
$A^{\phi}_*\egal A^\phi(f^{*}_\phi)$ the minimal $\phi-$risk. In
most of the cases studied $f^{*}_\phi$ or its sign is equal to the
Bayes classifier
$$ f^*(x)= {\rm{sign}}(2\eta(x)-1),$$ where $\eta$ is the
conditional probability function $x\longmapsto\mathbb{P}(Y=1|X=x)$
defined on $\cX$ (cf. \cite{bjm:06,lv:04,z:04}). The Bayes
classifier $f^*$ is a minimizer of the $\phi_0-$risk (cf.
\cite{dgl:96}).

Our  framework is the same as the one considered, among others, by
\cite{n:00,yang:00,catbook:01} and \cite{tsy:03,jrt:06}. We have a
family $\cF$ of $M$ classifiers $f_1,\ldots,f_M$ and a loss function
$\phi$. Our goal is to mimic the oracle
$\min_{f\in\cF}(A^{\phi}(f)-A^{\phi}_*)$ based on a sample $D_n$ of
$n$ i.i.d. observations $(X_1,Y_1),\ldots,(X_n,Y_n)$ of $(X,Y)$.
These classifiers may have been constructed from a previous sample
or they can belong to a dictionary of simple prediction rules like
decision stumps. The problem is to find a strategy which mimics as
fast as possible the best classifier in $\cF$. Such strategies can
then be used to construct efficient adaptive estimators (cf.
\cite{n:00,lec:05,lec4:06,cl:06}).  We consider the following
definition, which is inspired by the one given in \cite{tsy:03} for
the regression model.
\begin{definition}\label{SubERMdefoptimality}
  Let $\phi$ be a loss function. The
  remainder term $\gamma(n,M)$ is called {\bf optimal rate
  of aggregation for the $\phi-$risk},
  if the following two inequalities hold.
  \begin{enumerate}[i)]
         \item For any finite set $\cF$ of $M$ functions from $\cX$ to $[-1,1]$, there exists a statistic
          $\tilde{f}_n$ such that for any underlying probability measure $\pi$ and any integer $
            n\geq 1,$
            \begin{equation}\label{SubERMDefOracle}\mathbb{E}[A^{\phi}(\tilde{f}_n)-A^{\phi}_*]\leq
            \min_{f\in\cF}\left(A^{\phi}(f)-A^{\phi}_*
            \right)+C_1\gamma(n,M).\end{equation}
         \item There exists a finite set $\cF$ of $M$ functions from $\cX$ to $[-1,1]$
          such that for any statistic $\bar{f}_n$ there exists a
          probability distribution $\pi$ such that for all $n\geq1$
            \begin{equation}\label{SubERMDefLower}\mathbb{E}\left[A^{\phi}(\bar{f}_n)-A^{\phi}_* \right]\geq
            \min_{f\in\cF}\left(A^{\phi}(f)-A^{\phi}_*
            \right)+C_2\gamma(n,M).\end{equation}
  \end{enumerate}Here $C_1$ and $C_2$ are absolute positive constants which may depend on $\phi$.
   Moreover, when the above two properties i) and ii) are satisfied,
   we say that the procedure $\tilde{f}_n$, appearing in (\ref{SubERMDefOracle}), is an {\bf optimal
   aggregation procedure for the $\phi-$risk}.
\end{definition}

The paper is organized as follows. In the next Section we present
three aggregation strategies that will be shown to attain the
optimal rates of aggregation. Section $3$ presents performance of
these procedures. In Section $4$ we give some proofs of the
optimality of these procedures depending on the loss function. In
Section $5$ we state a result on suboptimality of the penalized
Empirical Risk Minimization procedures and of procedures called
selectors. In Section $6$ we give some remarks. All the proofs are
postponed to the last Section.

\section{Aggregation Procedures}\label{SubERMSectionAggregationProcedrues}
We introduce procedures that will be shown to achieve optimal rates
of aggregation depending on the loss function
$\phi:\mathbb{R}\longmapsto\mathbb{R}$. All these procedures are
constructed with the empirical version of the $\phi-$risk and the
main idea is that a classifier $f_j$ with a small empirical
$\phi-$risk is likely to have a small $\phi-$risk. We denote by
$$A_n^{\phi}(f)=\frac{1}{n}\sum_{i=1}^{n}\phi(Y_if(X_i))$$ the
empirical $\phi-$risk of a real-valued classifier $f$.

  The Empirical Risk
  Minimization ({\bf ERM}) procedure, is defined by  \begin{equation}\label{SubERMERM}
  \tilde{f}_n^{ERM}\in {\rm
  Arg}\min_{f\in\cF}A_n^{\phi}(f).\end{equation}This is an
  example of what we call a {\bf selector} which is an aggregate
  with values in the family $\cF$. Penalized ERM procedures are also
  examples of selectors.

  The Aggregation with
  Exponential Weights ({\bf AEW}) procedure  is given by
 \begin{equation}\label{SubERMAEW}\tilde{f}_n^{AEW}=\sum_{f\in\cF}w^{(n)}(f)f,\end{equation}
  where the
  weights $w^{(n)}(f)$ are defined by
  \begin{equation}\label{SubERMweightsAEW}w^{(n)}(f)=\frac{\exp\left( -nA_n^{\phi}(f)\right)}{\sum_{g\in\cF}\exp\left(
 -nA_n^{\phi}(g)\right)},\quad \forall f\in\cF.\end{equation}

 The Cumulative Aggregation with Exponential
 Weights ({\bf CAEW}) procedure, is defined by
 \begin{equation}\label{SubERMCAEW}\tilde{f}_{n,\beta}^{CAEW}=
 \frac{1}{n}\sum_{k=1}^n\tilde{f}_{k,\beta}^{AEW},\end{equation}
 where $\tilde{f}_{k,\beta}^{AEW}$ is constructed as in (\ref{SubERMAEW})
 based on
 the sample $(X_1,Y_1),\ldots,(X_k,Y_k)$ of size $k$ and with the 'temperature' parameter
 $\beta>0$. Namely,
 $$\tilde{f}_{k,\beta}^{AEW}=\sum_{f\in\cF}w^{(k)}_\beta(f)f, \mbox{ where }
 w^{(k)}_\beta(f)=\frac{\exp\left( -\beta^{-1}kA_k^{\phi}(f)\right)}{\sum_{g\in\cF}\exp\left(
 -\beta^{-1}kA_k^{\phi}(g)\right)},\quad \forall f\in\cF.$$

The idea of the ERM procedure goes to Le Cam and Vapnik. Exponential
weights have been discussed, for example, in
\cite{abb:97,hkw:98,kw:99,yang:00,catbook:01,lb:04,z:00,a:06} or in
\cite{v:90,lcb:06} in the on-line prediction setup.

\par

\section{Exact Oracle Inequalities.}
We now recall some known upper bounds on the excess risk. The first
point of the following Theorem goes to \cite{vc:81}, the second
point can be found in \cite{jntv:05} or \cite{cl:06} and the last
point, dealing with the case of a $\beta-$convex loss function, is
Corollary 4.4 of \cite{jrt:06}.

\begin{theorem}\label{SubERMtheoremoracleclassification}
Let $\phi:\mathbb{R}\longmapsto\mathbb{R}$ be a bounded loss
function. Let $\cF$ be a family of $M$  functions $f_1,\ldots,f_M$
with values in $[-1,1]$, where $M\geq2$ is an integer.
\begin{enumerate}[i)]
    \item The Empirical Risk Minimization procedure $\tilde{f}_n=\tilde{f}_n^{ERM}$
        satisfies\begin{equation}\label{SubERMInegOracleERM}
        \mathbb{E}[A^\phi(\tilde{f}_n)-A^{\phi}_*]\leq
        \min_{f\in\cF}(A^\phi(f)-A^{\phi}_*)+ C\sqrt{\frac{\log
        M}{n}},\end{equation}where $C>0$ is a constant depending only on
        $\phi$.
    \item If $\phi$ is convex, then the CAEW procedure
        $\tilde{f}_n=\tilde{f}_n^{CAEW}$ with ``temperature
        parameter''
        $\beta=1$ and the AEW procedure
        $\tilde{f}_n=\tilde{f}_n^{AEW}$ satisfy
        (\ref{SubERMInegOracleERM}).
    \item If $\phi$ is $\beta-$convex for a positive number $\beta$, then
        the CAEW procedure with ``temperature
        parameter'' $\beta$, satisfies
        $$\mathbb{E}[A^{\phi}(\tilde{f}_{n,\beta}^{CAEW})-A^{\phi}_*]\leq
        \min_{f\in\cF}(A^{\phi}(f)-A^{\phi}_*)+\beta\frac{\log M}{n}.$$
\end{enumerate}

\end{theorem}

\section{Optimal Rates of Aggregation.}
To understand how behaves the optimal rate of aggregation depending
on the loss we introduce a ``continuous scale'' of loss functions
indexed by a non negative number $h$,
\begin{equation*}
\phi_h(x)=\left\{
  \begin{array}{ll}    h\phi_1(x)+(1-h)\phi_0(x)      & \mbox{if } 0\leq h\leq 1\\
                       (h-1)x^2-x+1                  &     \mbox{if } h>1,\\
\end{array} \right.
\end{equation*}defined for any $x\in\mathbb{R}$, where $\phi_0$ is the $0-1$ loss
and $\phi_1$ is the hinge loss.


This set of losses is representative enough since it describes
different type of convexity: for any $h>1$, $\phi_h$ is
$\beta-$convex on $[-1,1]$ with $\beta\geq \beta_h\egal
(2h-1)^2/(2(h-1))\geq2$, for $h=1$ the loss is linear and for $h<1$,
$\phi_h$ is non-convex. For $h\geq0$, we consider
$$A_h(f)\egal A^{\phi_h}(f), f_h^*\egal f_{\phi_h}^* \mbox{ and } A_h^*\egal A^{\phi_h}_*=A^{\phi_h}(f_h^*).$$

\begin{theorem}\label{SubERMTheoLowerBounds}
Let $M\geq2$ be an integer. Assume that the space $\cX$ is infinite.

If $0\leq h<1$, then the optimal rate of aggregation for the
$\phi_h-$risk is achieved by the ERM procedure and is equal to
$$\sqrt{\frac{\log M}{n}}.$$

For $h=1$, the optimal rate of aggregation for the $\phi_1-$risk is
achieved by the ERM, the AEW and the CAEW (with 'temperature'
parameter $\beta=1$) procedures and is equal to
$$\sqrt{\frac{\log M}{n}}.$$

If $h>1$ then, the optimal rate of aggregation for the $\phi_h-$risk
is achieved by the CAEW, with 'temperature' parameter $\beta_h$ and
is equal to$$\frac{\log M}{n}.$$
\end{theorem}

\section{Suboptimality of Penalized ERM Procedures.}

In this Section we prove a lower bound under the margin assumption
for any selector and we give a more precise lower bound for
penalized ERM procedures. First, we recall the definition of the
margin assumption introduced in \cite{tsy:04}.

\noindent{\bf Margin Assumption(MA):} {\it The probability measure
$\pi$ satisfies the margin assumption MA($\kappa$), where
$\kappa\geq1$ if we have
\begin{equation}\label{MarginAssumption}\mathbb{E}[|f(X)-f^*(X)|]\leq c(A_0(f)-A_0^*)^{1/\kappa},\end{equation}for any measurable function $f$ with values in
$\{-1,1\}$}

\noindent We denote by $\cP_\kappa$ the set of all probability
distribution $\pi$ satisfying MA($\kappa$).

\begin{theorem}\label{SubERMTheoWeaknessERM}
Let $M\geq2$ be an integer, $\kappa\geq1$ be a real number, $\cX$ be
infinite and $\phi:\mathbb{R}\longmapsto\mathbb{R}$ be a loss
function such that $a_\phi\egal\phi(-1)-\phi(1)>0$. There exists a
family $\cF$ of $M$ classifiers with values in $\{-1,1\}$ satisfying
the following.

Let $\tilde{f}_n$ be a selector with values in  $\cF$. Assume that
$\sqrt{(\log M)/n}\leq1/2$. There exists a probability measure
$\pi\in\cP_\kappa$ and an absolute constant $C_3>0$ such that
$\tilde{f}_n$ satisfies
\begin{equation}\label{LBSelectors}\mathbb{E}\left[A^{\phi}(\tilde{f}_n)-A^{\phi}_* \right]\geq
 \min_{f\in\cF}\left(A^{\phi}(f)-A^{\phi}_*
 \right)+C_3\Big(\frac{\log
 M}{n}\Big)^{\frac{\kappa}{2\kappa-1}}.\end{equation}

Consider the penalized ERM procedure $\tilde{f}_n^{pERM}$ associated
with $\cF$, defined by
$$\tilde{f}_n^{pERM}\in{\rm Arg}\min_{f\in\cF}(A^{\phi}_n(f)+{\rm
pen}(f))$$ where the penalty function ${\rm{pen}}(\cdot)$ satisfies
$|{\rm pen}(f)|\leq C\sqrt{(\log M)/n},\forall f\in\cF,$ with $0\leq
C<\sqrt{2}/3$. Assume that $1188\pi C^2 M^{9C^2}\log M\leq n$. If
$\kappa>1$ then, there exists a probability measure
$\pi\in\cP_\kappa$ and an absolute constant $C_4>0$ such that the
penalized ERM procedure $\tilde{f}_n^{pERM}$ satisfies
$$\mathbb{E}\left[A^{\phi}(\tilde{f}_n^{pERM})-A^{\phi}_* \right]\geq
\min_{f\in\cF}\left(A^{\phi}(f)-A^{\phi}_*
\right)+C_4\sqrt{\frac{\log M}{n}}.$$
\end{theorem}

\begin{rem}
Inspection of the proof shows that Theorem
\ref{SubERMTheoWeaknessERM} is valid for any family $\cF$ of
classifiers $f_1,\ldots,f_M$, with values in $\{-1,1\}$,  such that
there exist points $x_1,\ldots,x_{2^M}$ in $\cX$ satisfying
$\big\{(f_1(x_j),\ldots,f_M(x_j)):j=1,\ldots,2^M\big\}=\{-1,1\}^M$.
\end{rem}

\begin{rem}
If we use a penalty function such that $|{\rm pen}(f)|\leq\gamma
n^{-1/2},\forall f\in\cF$, where $\gamma>0$ is an absolute constant
(i.e. $0\leq C \leq \gamma(\log M)^{-1/2}$), then the condition
``$1188\pi C^2 M^{9C^2}\log M\leq n$'' of Theorem
\ref{SubERMTheoWeaknessERM} is equivalent to ``$n$ greater than a
constant''.
\end{rem}

Theorem \ref{SubERMTheoWeaknessERM} states that the ERM procedure
(and even penalized ERM procedures) cannot mimic the best classifier
in $\cF$ with rates faster than $((\log M)/n)^{1/2}$ if the basis
classifiers in $\cF$ are different enough, under a very mild
condition on the loss $\phi$. If there is no margin assumption
(which corresponds to the case $\kappa=+\infty$), the result of
Theorem \ref{SubERMTheoWeaknessERM} can be easily deduced from the
lower bound in Chapter 7 of \cite{dgl:96}. The main message of
Theorem \ref{SubERMTheoWeaknessERM} is that such a negative
statement remains true even under the margin assumption
MA($\kappa$). Selectors aggregate cannot mimic the oracle faster
than $((\log M)/n)^{1/2}$ in general. Under MA($\kappa$), they
cannot mimic the best classifier in $\cF$ with rates faster than
$((\log M)/n)^{\kappa/(2\kappa-1)}$ (which is greater than $(\log
M)/n$ when $\kappa>1$). We know, according to Theorem
\ref{SubERMtheoremoracleclassification}, that the CAEW procedure
mimics the best classifier in $\cF$ at the rate $(\log M)/n$ if the
loss is $\beta-$convex. Thus, penalized ERM procedures (and more
generally, selectors) are suboptimal aggregation procedures when the
loss function is $\beta-$convex even if we add the constraint that
$\pi$ satisfies MA($\kappa$).

We can extend Theorem \ref{SubERMTheoWeaknessERM} to a more general
framework \cite{lec6:06} and we obtain that, if the loss function
associated with a risk is somewhat more than convex then it is
better to use aggregation procedures with exponential weights
instead of selectors (in particular penalized ERM or pure ERM). We
do not know whether the lower bound (\ref{LBSelectors}) is sharp,
i.e., whether there exists a selector attaining the reverse
inequality with the same rate.
\par

\section{Discussion.}
We proved in Theorem \ref{SubERMTheoLowerBounds} that the ERM
procedure is optimal only for non-convex losses and for the
borderline case of the hinge loss. But, for non-convex losses, the
implementation of the ERM procedure requires  minimization of a
function which is not convex. This is hard to implement and not
efficient from a practical point of view. In conclusion, the ERM
procedure is theoretically optimal only for non-convex losses but in
that case it is practically inefficient and it is practically
efficient only for the cases where ERM is theoretically suboptimal.

For any convex loss $\phi$, we have
$\frac{1}{n}\sum_{k=1}^nA^{\phi}(\tilde{f}^{AEW}_{k,\beta})\leq
A^{\phi}({\tilde{f}^{CAEW}_\beta})$. Next, less observations are
used for the construction of $\tilde{f}^{AEW}_{k,\beta},1\leq k \leq
n-1,$ than for the construction of $\tilde{f}^{AEW}_{n,\beta}$. We
can therefore expect the $\phi-$risk of $\tilde{f}^{AEW}_{n,\beta}$
to be smaller than the $\phi-$risk of $\tilde{f}^{AEW}_{k,\beta}$
for all $1\leq k\leq n-1$ and hence smaller than the $\phi-$risk of
$\tilde{f}^{CAEW}_{n,\beta}$.  Thus, the AEW procedure is likely to
be an optimal aggregation procedure for the convex loss functions.

The hinge loss happens to be really hinge for different reasons. For
losses "between" the $0-1$ loss and the hinge loss ($0\leq h\leq1$),
the ERM is an optimal aggregation procedure and the optimal rate of
aggregation is $\sqrt{(\log M)/n}$. For losses "over" the hinge loss
($h>1$), the ERM  procedure is suboptimal and $(\log M)/n$ is the
optimal rate of aggregation. Thus, there is a breakdown point in the
optimal rate of aggregation just after the hinge loss. This
breakdown can be explained by the concept of margin : this argument
has not been introduced here by the lack of space, but  can be found
in \cite{lec6:06}. 
Moreover for the hinge loss we get, by linearity $$\min_{f\in\cC}
A_1(f)-A_1^*=\min_{f\in\cF}A_1(f)-A^*_1,$$where $\cC$ is the convex
hull of $\cF$. Thus, for the particular case of the hinge loss,
``model selection'' aggregation and ``convex'' aggregation are
identical problems (cf. \cite{lec5:05} for more details).

\section{Proofs.}
{\bf{Proof of Theorem \ref{SubERMTheoLowerBounds}:}} The optimal
rates of aggregation of Theorem \ref{SubERMTheoLowerBounds} are
achieved by the procedures introduced in Section
\ref{SubERMSectionAggregationProcedrues}. Depending on the value of
$h$, Theorem \ref{SubERMtheoremoracleclassification} provides the
exact oracle inequalities required by the point
(\ref{SubERMDefOracle}) of Definition \ref{SubERMdefoptimality}. To
show optimality of these rates of aggregation, we need only to prove
the corresponding lower bounds. We consider two cases: $0\leq
h\leq1$ and $h>1$. Denote by $\cP$ the set of all probability
distributions on $\cX\times\{-1,1\}$.

Let $0\leq h \leq1$. It is easy to check that the Bayes rule $f^*$
is a minimizer of the $\phi_h-$risk. Moreover, using the inequality
$A_1(f)-A_1^*\geq A_0(f)-A_0^*,$ which holds for any real-valued
function $f$ (cf. \cite{z:04}), we have for any prediction rules
$f_1,\ldots,f_M$ (with values in $\{-1,1\}$) and for any finite set
$\cF$ of $M$ real valued functions,
\begin{eqnarray}\label{SubERMIneqBornInf}
\lefteqn{\inf_{\hat{f}_n}\sup_{\pi\in\cP}
\left(\mathbb{E}\left[A_h(\hat{f}_n)-A^*_h
\right]-\min_{f\in\cF}(A_h(f)-A^*_h) \right)}\\ & \geq &
\inf_{\hat{f}_n}\sup_{\substack{\pi\in\cP\\
f^*\in\{f_1,\ldots,f_M\}}} \Big(\mathbb{E}\left[A_h(\hat{f}_n)-A^*_h
\right]\Big)  \geq
\inf_{\hat{f}_n}\sup_{\substack{\pi\in\cP\\
f^*\in\{f_1,\ldots,f_M\}}} \Big(\mathbb{E}\left[A_0(\hat{f}_n)-A^*_0
\right]\Big)\nonumber.\end{eqnarray}

Let $N$ be an integer such that $2^{N-1}\leq M$, $x_1,\ldots,x_N$ be
$N$ distinct points of $\cX$ and  $w$ be a positive number
satisfying $(N-1)w\leq1$. Denote by $P^X$ the probability measure on
$\cX$ such that $P^X(\{x_j\})=w$, for $j=1,\ldots,N-1$ and
$P^X(\{x_N\})=1-(N-1)w$. We consider the cube
$\Omega=\{-1,1\}^{N-1}$. Let $0<\mathfrak{h}<1$. For all
$\sigma=(\sigma_1,\ldots,\sigma_{N-1})\in\Omega$ we consider
$$\eta_{\sigma}(x)=\left\{\begin{array}{ll} (1+\sigma_j \mathfrak{h})/2 & \mbox{if } x=x_1,\ldots,x_{N-1},\\
1 & \mbox{if } x =x_N. \end{array}
 \right.$$
For all $\sigma\in\Omega$ we denote by $\pi_\sigma$ the probability
measure on $\cX\times\{-1,1\}$ defined by its marginal $P^X$ on
$\cX$ and its conditional probability function $\eta_{\sigma}$.

We denote by $\rho$ the Hamming distance on $\Omega$. Let $\sigma,
\sigma'\in\Omega$ such that $\rho(\sigma,\sigma')=1$. Denote by $H$
the Hellinger's distance. Since $H^2\left(\pi_\sigma^{\otimes n},
\pi_{\sigma'}^{\otimes
n}\right)=2\Big(1-\Big(1-H^2(\pi_{\sigma},\pi_{\sigma'})/2\Big)^n\Big)$
and $H^2(\pi_{\sigma},\pi_{\sigma'})=2w(1-\sqrt{1-\mathfrak{h}^2}),$
then, the Hellinger's distance between the measures
$\pi_\sigma^{\otimes n}$ and $\pi_{\sigma'}^{\otimes n}$ satisfies
$$H^2\left(\pi_\sigma^{\otimes n}, \pi_{\sigma'}^{\otimes
n}\right)=2\left(1-(1-w(1-\sqrt{1-\mathfrak{h}^2}))^n\right).$$

Take $w$ and $\mathfrak{h}$ such that
$w(1-\sqrt{1-\mathfrak{h}^2})\leq n^{-1}.$ Then,
$H^2\left(\pi_\sigma^{\otimes n}, \pi_{\sigma'}^{\otimes
n}\right)\leq 2(1-e^{-1})<2$ for any integer $n$.

Let $\sigma\in\Omega$ and $\hat{f}_n$ be an estimator with values in
$\{-1,1\}$ (only the sign of a statistic is used when we work with
the $0-1$ loss). For $\pi=\pi_\sigma$, we have
$$
\mathbb{E}_{\pi_\sigma}[A_0(\hat{f}_n)-A_0^*]  \geq  \mathfrak{h}w
\mathbb{E}_{\pi_\sigma}\Big[\sum_{j=1}^{N-1}|\hat{f}_n(x_j)-\sigma_j|\Big].
$$Using Assouad's Lemma (cf. Lemma \ref{SubERMLemAssouadModif}), we obtain
\begin{equation}\label{SubERMIneqAssouadApplique}\inf_{\hat{f}_n}\sup_{\sigma\in\Omega}
\left(\mathbb{E}_{\pi_\sigma}\left[A_0(\hat{f}_n)-A_0^*
\right]\right)\geq \mathfrak{h}w\frac{N-1}{4e^2}.\end{equation}

Take now $w=(n\mathfrak{h}^2)^{-1}$, $N=\lceil \log M/ \log 2
\rceil$, $\mathfrak{h}=\left(n^{-1}\lceil \log M/ \log 2 \rceil
\right)^{1/2}$. We complete the proof by replacing $w$,
$\mathfrak{h}$ and $N$ in (\ref{SubERMIneqAssouadApplique}) and
(\ref{SubERMIneqBornInf}) by their values.

For the case $h>1$, we consider  an integer $N$ such that
$2^{N-1}\leq M$, $N-1$ different points $x_1,\ldots,x_N$ of $\cX$
and a positive number $w$ such that $(N-1)w\leq1$. We denote by
$P^X$ the probability measure on $\cX$ such that $P^X(\{x_j\})=w$
for $j=1,\ldots,N-1$ and $P^X(\{x_N\})=1-(N-1)w$. Denote by $\Omega$
the cube $\{-1,1\}^{N-1}$. For any $\sigma\in\Omega$ and $h>1$, we
consider the conditional probability function $\eta_\sigma$ in two
different cases. If $2(h-1)\leq 1$ we take
$$\eta_{\sigma}(x)=\left\{\begin{array}{ll} (1+2\sigma_j(h-1))/2 & \mbox{if }
 x=x_1,\ldots,x_{N-1}\\
2(h-1) & \mbox{if } x =x_N, \end{array}
 \right.$$and if $2(h-1)>1$ we take
$$\eta_{\sigma}(x)=\left\{\begin{array}{ll} (1+\sigma_j)/2 & \mbox{if }
 x=x_1,\ldots,x_{N-1}\\
1 & \mbox{if } x =x_N. \end{array}
 \right.$$For all $\sigma\in\Omega$ we denote by $\pi_\sigma$ the
probability measure on $\cX\times\{-1,1\}$ with the marginal $P^X$
on $\cX$ and the conditional probability function $\eta_{\sigma}$ of
$Y$ knowing $X$.

Consider $$\rho(h)=\left\{\begin{array}{cl} 1 & \mbox{if }
2(h-1)\leq1\\
(4(h-1))^{-1} & \mbox{if } 2(h-1)>1 \end{array}\right. \mbox{ and }
g^*_\sigma(x)=\left\{\begin{array}{cl}\sigma_j & \mbox{if } x=x_1,\ldots,x_{N-1}\\
1 & \mbox{if } x=x_N.\end{array} \right.$$ A minimizer of the
$\phi_h-$risk when the underlying distribution is $\pi_{\sigma}$ is
given by
$$f^*_{h,\sigma}\egal \frac{2\eta_{\sigma}(x)-1}{2(h-1)}=
\rho(h)g^*_\sigma(x),\quad \forall x\in \cX,$$ for any $h>1$ and
$\sigma\in\Omega.$

When we choose $\{f^*_{h,\sigma}:\sigma\in\Omega\}$ for the set
$\cF=\{f_1,\ldots,f_M\}$ of basis functions, we obtain
\begin{eqnarray*}
\sup_{\{f_1,\ldots,f_M\}}\inf_{\hat{f}_n}\sup_{\pi\in\cP}
\left(\mathbb{E}\left[A_h(\hat{f}_n)-A^*_h
\right]-\min_{j=1,\ldots,M}(A_h(f_j)-A^*_h) \right)\\ \geq
\inf_{\hat{f}_n}\sup_{\substack{\pi\in\cP:\\f^*_h\in\{f^*_{h,\sigma}:\sigma\in\Omega\}}}
\left(\mathbb{E}\left[A_h(\hat{f}_n)-A^*_h \right]
\right).\end{eqnarray*} Let $\sigma$ be an element of $\Omega$.
Under the probability distribution $\pi_\sigma$, we have
$A_h(f)-A^*_h=(h-1)\mathbb{E}[(f(X)-f^*_{h,\sigma}(X))^2],$ for any
real-valued function $f$ on $\cX$. Thus, for a real valued estimator
$\hat{f}_n$ based on $D_n$, we have
$$A_h(\hat{f}_n)-A^*_h\geq (h-1)w\sum_{j=1}^{N-1}(\hat{f}_n(x_j)-\rho(h)\sigma_j)^2.$$
We consider the projection function $\psi_h(x)=\psi(x/\rho(h))$ for
any $x\in\cX$, where $\psi(y)=\max(-1,\min(1,y)),\forall
y\in\mathbb{R}$. We have
\begin{eqnarray*}
\mathbb{E}_\sigma[A_h(\hat{f}_n)-A^*_h] & \geq & w(h-1)\sum_{j=1}^{N-1}\mathbb{E}_\sigma(\psi_h(\hat{f}_n(x_j))-\rho(h)\sigma_j)^2\\
   & \geq  & w(h-1)(\rho(h))^2\sum_{j=1}^{N-1}\mathbb{E}_\sigma(\psi(\hat{f}_n(x_j))-\sigma_j)^2\\
&  \geq&  4w(h-1)(\rho(h))^2
\inf_{\hat{\sigma}\in[0,1]^{N-1}}\max_{\sigma\in\Omega}\mathbb{E}_\sigma
\left[\sum_{j=1}^{N-1}\left|\hat{\sigma}_j-\sigma_j
\right|^2\right],
\end{eqnarray*}where the infimum
$\inf_{\hat{\sigma}\in[0,1]^{N-1}}$ is taken over all estimators
$\hat{\sigma}$ based on one observation from the statistical
experience $\left\{\pi_\sigma^{\otimes n}|\sigma\in\Omega \right\}$
and with values in $[0,1]^{N-1}$.

For any $\sigma,\sigma'\in\Omega$ such that
$\rho(\sigma,\sigma')=1,$ the Hellinger's distance between the
measures $\pi_\sigma^{\otimes n}$ and  $\pi_{\sigma'}^{\otimes n}$
satisfies
$$H^2\left(\pi_\sigma^{\otimes n}, \pi_{\sigma'}^{\otimes
n}\right)=\left\{\begin{array}{lc}
2\left(1-(1-2w(1-\sqrt{1-h^2}))^n\right) & \mbox{if }2(h-1)<1\\
2\left(1-(1-2w(1-\sqrt{3/4}))^n\right) & \mbox{if }2(h-1)\geq1\\
\end{array}\right..$$

We take $$w=\left\{\begin{array}{cl}
(2n(h-1)^2) &\mbox{if }2(h-1)<1\\
8n^{-1}& \mbox{if }2(h-1)\geq1.\\
\end{array}\right.$$Thus, we have for
any $\sigma,\sigma'\in\Omega$ such that $\rho(\sigma,\sigma')=1,$
$$H^2\left(\pi_\sigma^{\otimes n}, \pi_{\sigma'}^{\otimes
n}\right)\leq 2(1-e^{-1}).$$

To complete the proof we apply Lemma \ref{SubERMLemAssouadModif}
with $N=\lceil (\log M)/n\rceil$.

{\bf{Proof of Theorem \ref{SubERMTheoWeaknessERM}:}} Consider $\cF$
a family of classifiers $f_1,\ldots,f_M$, with values in $\{-1,1\}$,
such that there exist $2^M$ points $x_1,\ldots,x_{2^M}$ in $\cX$
satisfying
$\big\{(f_1(x_j),\ldots,f_M(x_j)):j=1,\ldots,2^M\big\}=\{-1,1\}^M\egal\cS_M$.

Consider the lexicographic order on $\cS_M$:
$$(-1,\ldots,-1)\preccurlyeq(-1,\ldots,-1,1)\preccurlyeq(-1,\ldots,-1,1,-1)\preccurlyeq\ldots\preccurlyeq(1,\ldots,1).$$
Take $j$ in $\{1,\ldots,2^M\}$ and denote by $x_j'$ the element in
$\{x_1,\ldots,x_{2^M}\}$ such that $(f_1(x_j'),\ldots,f_M(x_j'))$ is
the $j-$th element of $\cS_M$ for the lexicographic order. We denote
by $\varphi$ the bijection between $\cS_M$ and
$\{x_1,\ldots,x_{2^M}\}$ such that the value of $\varphi$ at the
$j-$th element of $\cS_M$ is $x_j'$. By using the bijection
$\varphi$ we can work independently either on the set $\cS_M$ or on
$\{x_1,\ldots,x_{2^M}\}$. Without any assumption on the space $\cX$,
we consider, in what follows, functions and probability measures on
$\cS_M$. Remark that for the bijection $\varphi$ we have
$$f_j(\varphi(x))=x^j,\quad \forall x=(x^1,\ldots,x^M)\in\cS_M, \forall j\in\{1,\ldots,M\}.$$
With a slight abuse of notation, we still denote by $\cF$ the set of
functions $f_1,\ldots,f_M$ defined by $f_j(x)=x^j,$ for any
$j=1,\ldots,M.$

First remark that for any $f,g$ from $\cX$ to $\{-1,1\}$, using $
 \mathbb{E}[\phi(Yf(X))|X]=\mathbb{E}[\phi(Y)|X]\1_{(f(X)=1)}+\mathbb{E}[\phi(-Y)|X]\1_{(f(X)=-1)},
$ we have $$\mathbb{E}[\phi(Y f(X))|X]-\mathbb{E}[\phi(Y
g(X))|X]=a_\phi(1/2-\eta(X))(f(X)-g(X)).$$ Hence, we obtain
$A^{\phi}(f)-A^{\phi}(g)=a_\phi(A_0(f)-A_0(g)).$ So, we have for any
$j=1,\ldots,M,$
$$A^{\phi}(f_j)-A^{\phi}(f^*)=a_\phi(A_0(f_j)-A_0^*).$$
Moreover, for any $f:\cS_M\longmapsto\{-1,1\}$ we have
$A_n^{\phi}(f)=\phi(1)+a_\phi A_n^{\phi_0}(f)$ and $a_\phi>0$ by
assumption, hence,
$$\tilde{f}^{pERM}_n\in{\rm{Arg}}\min_{f\in\cF}(A_n^{\phi_0}(f)+{\rm
pen}(f)).$$Thus, it suffices to prove Theorem
\ref{SubERMTheoWeaknessERM}, when the loss function $\phi$ is the
classical $0-1$ loss function $\phi_0$.

We denote by $\cS_{M+1}$ the set $\{-1,1\}^{M+1}$  and by
$X^{0},\ldots,X^{M}$, $M+1$ independent random variables with values
in $\{-1,1\}$ such that $X^0$ is distributed according to a
Bernoulli $\cB(w,1)$ with parameter $w$  (that is
$\mathbb{P}(X^0=1)=w$ and $\mathbb{P}(X^0=-1)=1-w$) and  the $M$
other variables $X^1,\ldots,X^M$ are distributed according to a
Bernoulli $\cB(1/2,1)$. The parameter $0\leq w \leq 1$ will be
chosen wisely in what follows.

For any $j\in\{1,\ldots,M\}$, we consider the probability
distribution $\pi_j=(P^X,\eta^{(j)})$ of a couple of random
variables $(X,Y)$ with values in $S_{M+1}\times\{-1,1\}$, where
$P^X$ is the probability distribution on $S_{M+1}$ of
$X=(X^{0},\ldots,X^{M})$ and $\eta^{(j)}(x)$ is the regression
function at the point $x\in\cS_{M+1}$, of $Y=1$ knowing that $X=x$,
given by
$$\eta^{(j)}(x)=\left\{\begin{array}{ll}
1       & \mbox{if } x^0=1\\
1/2+h/2 & \mbox{if } x^0=-1, x^j=-1\\
1/2+h   & \mbox{if } x^0=-1, x^j=1\\
\end{array}\right., \quad \forall x=(x^0,x^1,\ldots,x^M)\in\cS_{M+1},$$
where $h>0$ is a parameter chosen wisely in what follows. The Bayes
rule $f^*$, associated with the distribution
$\pi_j=(P^X,\eta^{(j)})$, is identically equal to $1$ on
$\cS_{M+1}$.

If the probability distribution of $(X,Y)$ is $\pi_j$ for a
$j\in\{1,\ldots,M\}$ then, for any $0<t<1$, we have
$\mathbb{P}[|2\eta(X)-1|\leq t]\leq (1-w)\1_{h\leq t}.$ Now, we take
$$1-w=h^{\frac{1}{\kappa-1}},$$ then, we have
$\mathbb{P}[|2\eta(X)-1|\leq t]\leq t^{\frac{1}{\kappa-1}}$ and so
$\pi_j\in\cP_\kappa$.

We extend the definition of the $f_j$'s to the set $\cS_{M+1}$ by
$f_j(x)=x^j$ for any $x=(x^0,\ldots,x^M)\in \cS_{M+1}$ and
$j=1,\ldots,M$. Consider $\cF=\{f_1,\ldots,f_M\}.$ Assume that
$(X,Y)$ is distributed according to $\pi_j$ for a
$j\in\{1,\ldots,M\}$.  For any $k\in\{1,\ldots,M\}$ and $k\neq j$,
we have
\begin{equation*}
  A_0(f_k)-A_0^* = \sum_{x\in\cS_{M+1}}|\eta(x)-1/2||f_k(x)-1|\mathbb{P}[X=x]
  = \frac{3h(1-w)}{8}+\frac{w}{2}
\end{equation*}
and the excess risk of $f_j$ is given by
$A_0(f_j)-A_0^*=(1-w)h/4+w/2.$ Thus, we have
$$\min_{f\in\cF}A_0(f)-A_0^*=A_0(f_j)-A_0^*=(1-w)h/4+w/2.$$

First, we prove the lower bound for any selector. Let $\tilde{f}_n$
be a selector with values in $\cF$. If the underlying probability
measure is $\pi_j$ for a $j\in\{1,\ldots,M\}$ then,
\begin{align*}
  \mathbb{E}^{(j)}_n[A_0(\tilde{f}_n)-A_0^*]&=\sum_{k=1}^M(A_0(f_k)-A_0^*)\pi_j^{\otimes n}[\tilde{f}_n=f_k]\\
  &=\min_{f\in\cF}(A_0(f)-A_0^*)+
  \frac{h(1-w)}{8}\pi_j^{\otimes n}[\tilde{f}_n\neq f_j],
\end{align*}where $\mathbb{E}^{(j)}_n$ denotes the expectation
w.r.t. the observations $D_n$ when $(X,Y)$ is distributed according
to $\pi_j$. Hence, we have
\begin{equation*}
  \max_{1\leq j \leq
  M}\{\mathbb{E}^{(j)}_n[A_0(\tilde{f}_n)-A_0^*]-\min_{f\in\cF}(A_0(f)-A_0^*)\}\geq\frac{h(1-w)}{8}\inf_{\hat{\phi}_n}\max_{1\leq j \leq
  M}\pi_j^{\otimes n}[\hat{\phi}_n\neq j],
\end{equation*}
where the infimum $\inf_{\hat{\phi}_n}$ is taken over all tests valued in $\{1,\ldots,M\}$ constructed from one observation
in the model
$(\cS_{M+1}\times\{-1,1\},\cA\times\cT,\{\pi_1,\ldots,\pi_M\})^{\otimes
n},$ where $\cT$ is the natural $\sigma-$algebra on $\{-1,1\}$.
Moreover, for any $j\in\{1,\ldots,M\}$, we have
$$K(\pi_j^{\otimes n}|\pi_1^{\otimes n})\leq
\frac{nh^2}{4(1-h-2h^2)},$$where $K(P|Q)$ is the Kullback-Leibler
divergence between $P$ and $Q$ (that is $\int\log(dP/dQ)dP$ if
$P<<Q$ and $+\infty$ otherwise). Thus, if we apply Lemma
\ref{LowerBoundTest} with $h=((\log M)/n)^{(\kappa-1)/(2\kappa-1)}$,
we obtain the result.

Second, we prove the lower bound for the pERM procedure
$\hat{f}_n=\tilde{f}_n^{pERM}$. Now, we assume that the probability
distribution of $(X,Y)$ is $\pi_M$ and we take
\begin{equation}\label{EquaValueh}h=\Big(C^2\frac{\log M}{n}\Big)^{\frac{\kappa-1}{2\kappa}}.\end{equation}
We have $\DS
  \mathbb{E}[A_0(\hat{f}_n)-A_0^*]=\min_{f\in\cF}(A_0(f)-A_0^*)+
  \frac{h(1-w)}{8}\mathbb{P}[\hat{f}_n\neq f_M].
$ Now, we upper bound $\mathbb{P}[\hat{f}_n= f_M]$, conditionally to
$\cY=(Y_1,\ldots,Y_n)$. We have
\begin{eqnarray*}
  \lefteqn{\mathbb{P}[\hat{f}_n= f_M|\cY]}\\&= &\mathbb{P}[\forall j=1,\ldots,M-1, A_n^{\phi_0}(f_M)+
  {\rm{pen}}(f_M)\leq A_n^{\phi_0}(f_j)+{\rm{pen}}(f_j)|\cY]\\
  &=&\mathbb{P}[\forall j=1,\ldots,M-1, \nu_M
  \leq \nu_j+n({\rm{pen}}(f_j)-{\rm{pen}}(f_M))|\cY],
\end{eqnarray*} where  $\nu_j=\sum_{i=1}^n\1_{(Y_iX_{i}^j\leq0)},\forall j=1,\ldots,M$ and
$X_i=(X_{i}^j)_{j=0,\ldots,M}\in\cS_{M+1},\forall i=1,\ldots,n$.
Moreover, the coordinates $X_{i}^j,i=1,\ldots,n;j=0,\ldots,M$ are
independent, $Y_1,\ldots,Y_n$ are independent of
$X_{i}^j,i=1,\ldots,n;j=1,\ldots,M-1$ and $|{\rm{pen}}(f_j)|\leq
h^{\kappa/(\kappa-1)},\forall j=1,\ldots,M$. So, we have
\begin{eqnarray*}
 \mathbb{P}[\hat{f}_n= f_M|\cY]&=&\sum_{k=0}^n\mathbb{P}[\nu_M=k|\cY]\prod_{j=1}^{M-1}\mathbb{P}[k\leq
  \nu_j+n({\rm{pen}}(f_j)-{\rm{pen}}(f_M))|\cY]
   \\
    &\leq & \sum_{k=0}^n\mathbb{P}[\nu_M=k|\cY]\Big(\mathbb{P}[k\leq
  \nu_1+2nh^{\kappa/(\kappa-1)}|\cY]\Big)^{M-1}
  \\
  &\leq&\mathbb{P}[\nu_M\leq \bar{k}|\cY]+\big(\mathbb{P}[\bar{k}\leq
  \nu_1+2nh^{\kappa/(\kappa-1)}|\cY]\big)^{M-1},
\end{eqnarray*}where \begin{eqnarray*}
\bar{k}&=&\mathbb{E}[\nu_M|\cY]-2nh^{\kappa/(\kappa-1)}\\
       &=&\frac{1}{2}\sum_{i=1}^n\Big(\frac{2-4h}{2-3h}\1_{(Y_i=-1)}
+\frac{1+h^{1/(\kappa-1)}(h/2-1/2)}{1+h^{1/(\kappa-1)}(3h/4-1/2)}\1_{(Y_i=1)}\Big)-2nh^{\kappa/(\kappa-1)}.
\end{eqnarray*}
Using Einmahl and Masson's concentration inequality (cf.
\cite{em:96}), we obtain
$$\mathbb{P}[\nu_M\leq \bar{k}|\cY]\leq
\exp(-2nh^{2\kappa/(\kappa-1)}).$$ Using Berry-Ess{\'e}en's theorem
(cf. p.471 in  \cite{bd:01}), the fact that $\cY$ is independent of
$(X_{i}^j;1\leq i \leq n, 1\leq j\leq M-1)$ and $\bar{k}\geq
n/2-9nh^{\kappa/(\kappa-1)}/4$, we get
\begin{equation*}
\mathbb{P}[\bar{k}\leq
  \nu_1+2nh^{\frac{\kappa}{\kappa-1}}|\cY]\leq \mathbb{P}\left[\frac{n/2-\nu_1}{\sqrt{n}/2}\leq
  6h^{\frac{\kappa}{\kappa-1}}\sqrt{n}\right]
  \leq\Phi(6h^{\frac{\kappa}{\kappa-1}}\sqrt{n})+\frac{66}{\sqrt{n}},
\end{equation*}where ${\Phi}$ stands for the standard normal
distribution function. Thus, we have
\begin{eqnarray}\label{EquaLowerBoundERMFinal1}
  \lefteqn{\mathbb{E}[A_0(\hat{f}_n)-A_0^*] \geq
  \min_{f\in\cF}(A_0(f)-A_0^*)}\\
  &&+\frac{(1-w)h}{8}\Big(1-\exp(-2nh^{2\kappa/(\kappa-1)})-\Big(\Phi(6h^{\kappa/(\kappa-1)}\sqrt{n})
  +66/\sqrt{n}\Big)^{M-1} \Big)\nonumber.
\end{eqnarray}

Next, for any $a>0$, by the elementary properties of the tails of
normal distribution, we have
\begin{equation}\label{LemQueueGaussienne}
  1-\Phi(a)=\frac{1}{\sqrt{2\pi}}\int_{a}^{+\infty}\exp(-t^2/2)dt\geq\frac{a}{\sqrt{2\pi}(a^2+1)}e^{-a^2/2}.
\end{equation}
Besides, we have for $0<C<\sqrt{2}/6$ (a modification for $C=0$ is
obvious) and $(3376C)^2(2\pi M^{36C^2}\log M)\leq n$, thus, if we
replace $h$ by its value given in (\ref{EquaValueh}) and if we apply
(\ref{LemQueueGaussienne}) with $a=16C\sqrt{\log M}$, then we obtain
\begin{equation}\label{EquaLowerBoundERMFinal2}
  \Big(\Phi(6h^{\kappa/(\kappa-1)}\sqrt{n})+66/\sqrt{n}\Big)^{M-1}
\leq \exp\Big[-\frac{M^{1-18C^2}}{18C\sqrt{2\pi\log
M}}+\frac{66(M-1)}{\sqrt{n}}\Big].
\end{equation}
Combining (\ref{EquaLowerBoundERMFinal1}) and
(\ref{EquaLowerBoundERMFinal2}), we obtain the result with
$C_4=(C/4)\Big(1-\exp(-8C^2)-\exp(-1/(36C\sqrt{2\pi\log2}))\Big)>0.$
\begin{flushright}$\blacksquare$\end{flushright}

The following lemma is used to establish the lower bounds of Theorem
\ref{SubERMTheoLowerBounds}. It is a version of Assouad's Lemma (cf.
\cite{tsybook:04}). Proof can be found in \cite{lec6:06}.

\begin{lemma}\label{SubERMLemAssouadModif}
Let $(\cX,\cA)$ be a measurable space. Consider a set of probability
$\{P_{\omega}/\omega\in\Omega\}$ indexed by the cube
$\Omega=\{0,1\}^m$. Denote by $\mathbb{E}_\omega$ the expectation
under $P_\omega$. Let $\theta\geq1$ be a number. Assume that:
$$\forall \omega,\omega' \in\Omega/ \rho(\omega,\omega')=1,\
 H^2(P_{\omega},P_{\omega'})\leq \alpha <2,$$
then we have
$$\inf_{\hat{w}\in[0,1]^m}\max_{\omega\in\Omega}\mathbb{E}_\omega
\left[\sum_{j=1}^m\left|\hat{w_j}-w_j \right|^\theta\right]\geq
m2^{-3-\theta}(2-\alpha)^2$$ where the infimum
$\inf_{\hat{w}\in[0,1]^m}$ is taken over all estimator based on an
observation from the statistical experience
$\left\{P_\omega|\omega\in\Omega \right\}$ and with values in
$[0,1]^m$.
\end{lemma}

We use the following lemma to prove the weakness of selector
aggregates. A proof can be found p. 84 in \cite{tsybook:04}.

\begin{lemma}\label{LowerBoundTest}
Let $\mathbb{P}_1,\ldots,\mathbb{P}_M$ be $M$ probability measures
on a measurable space $(\cZ,\cT)$ satisfying
$\DS\frac{1}{M}\sum_{j=1}^M K(\mathbb{P}_j|\mathbb{P}_1)\leq \alpha
\log M,$ where $0<\alpha<1/8$. We have
$$\inf_{\hat{\phi}}\max_{1\leq j \leq M} \mathbb{P}_j(\hat{\phi}\neq
j)\geq
\frac{\sqrt{M}}{1+\sqrt{M}}\Big(1-2\alpha-2\sqrt{\frac{\alpha}{\log
2}}\Big),
$$ where the infimum $\inf_{\hat{\phi}}$ is taken over all tests $\hat{\phi}$ with
values in $\{1,\ldots,M\}$ constructed from one observation in the
statistical model $(\cZ,\cT,\{\mathbb{P}_1,\ldots,\mathbb{P}_M\}).$
\end{lemma}

\par

\end{document}